\documentclass[10pt]{article}
\usepackage{mathrsfs}
\usepackage{amssymb}
\usepackage{amsmath}
\usepackage[all]{xy}

\def\Im{\mathop{\rm Im}\nolimits}
\def\Ker{\mathop{\rm Ker}\nolimits}
\def\Coker{\mathop{\rm Coker}\nolimits}

\def\mod{\mathop{\rm mod}\nolimits}
\def\Mod{\mathop{\rm Mod}\nolimits}
\def\Hom{\mathop{\rm Hom}\nolimits}

\def\Ext{\mathop{\rm Ext}\nolimits}

\def\Ab{\mathop{\rm Ab}\nolimits}

\title{\large \bf Applications of Exact Structures in Abelian Categories
}
\author{Junfu Wang\thanks{{\it E-mail address}: wangjunfu05@126.com}, Huanhuan Li\thanks{{\it E-mail address}: lihuanhuan0416@163.com}
\ and Zhaoyong Huang\thanks{{\it E-mail address}: huangzy@nju.edu.cn}\\
{\footnotesize \it Department of Mathematics, Nanjing
University, Nanjing 210093, Jiangsu Province, China}}
\date{}
\begin{document}
\baselineskip=16pt \maketitle

\begin{abstract}  In an abelian category $\mathscr{A}$ with small $\Ext$ groups, we show that there exists a one-to-one correspondence between any two of the following:
balanced pairs, subfunctors $\mathcal{F}$ of $\Ext^{1}_{\mathscr{A}}(-,-)$ such that $\mathscr{A}$ has enough $\mathcal{F}$-projectives
and enough $\mathcal{F}$-injectives and Quillen exact structures $\mathcal{E}$ with enough $\mathcal{E}$-projectives and enough $\mathcal{E}$-injectives.
In this case, we get a strengthened version of the translation of the Wakamatsu lemma to the exact context, and also prove that subcategories
which are $\mathcal{E}$-resolving and epimorphic precovering with kernels in their right $\mathcal{E}$-orthogonal class and subcategories
which are $\mathcal{E}$-coresolving and monomorphic preenveloping with cokernels in their left $\mathcal{E}$-orthogonal class are determined by
each other. Then we apply these results to construct some (pre)enveloping and (pre)covering classes and complete hereditary $\mathcal{E}$-cotorsion pairs
in the module category.
\end{abstract}

\vspace{0.5cm}

{\bf  1. Introduction}

\vspace{0.2cm}

Throughout this paper, $\mathscr{A}$ is an abelian category and a subcategory of $\mathscr{A}$ means a full and additive subcategory closed under isomorphisms and direct summands.

The notion of exact categories was originally due to Quillen [17]. The theory of exact categories was developed by B\"{u}hler, Fu, Gillespie, Hovey, Keller, Krause, Neeman,
\v{S}\v{t}ov\'{i}\v{c}ek and possibly others, see [4, 8, 10, 13, 14, 15, 16, 20, 21], and so on. In addition, Auslander and Solberg developed in [1, 2] the theory of relative homological algebra with respect to a
subfunctor $\mathcal{F}$ of $\Ext^{1}_{\mathscr{A}}(-,-): \mathscr{A}^{op}\times \mathscr{A} \to \Ab$, where $\Ab$ is the category of abelian groups and $\mathscr{A}=\mod\mathrm{\Lambda}$
is the category of finitely generated modules over an artin algebra $\mathrm{\Lambda}$. Then, as a special case, Chen introduced and studied in [5] relative homology with respect to
balanced pairs in an abelian category. On the other hand, the notion of cotorsion pairs was introduced by Salce in [19], which is based on the functor $\Ext^{1}_{\mathscr{A}}(-,-)$.
The theory of cotorsion pairs, studied by many authors, has played an important role in homological algebra and representation theory of algebras, see [7, 8, 9, 11, 13, 15, 20, 21] and references therein.
Motivated by these, in this paper, for an abelian category $\mathscr{A}$ with small $\Ext$ groups, we first investigate the relations among exact structures,
subfunctors of $\Ext^{1}_{\mathscr{A}}(-,-)$ and balanced pairs in $\mathscr{A}$, then we study cotorsion pairs with respect to exact structures in $\mathscr{A}$. The paper is organized as follows.

In Section 2, we recall the definitions of exact categories, balanced pairs and some notions in relative homological algebra. We prove that there exists a one-to-one correspondence between
any two of the following: (1) Balanced pairs $(\mathscr{C},\mathscr{D})$ in $\mathscr{A}$;
(2) Subfunctors $\mathcal{F}\subseteq \Ext^{1}_{\mathscr{A}}(-,-)$ such that $\mathscr{A}$ has enough $\mathcal{F}$-projectives and enough $\mathcal{F}$-injectives;
(3) Quillen exact structures $\mathcal{E}$ in $\mathscr{A}$ with enough $\mathcal{E}$-projectives and enough $\mathcal{E}$-injectives.

Let $\mathcal{E}$ be an exact structure on $\mathscr{A}$ such that $\mathscr{A}$ has enough $\mathcal{E}$-projectives and enough $\mathcal{E}$-injectives. In Section 3,
we get a strengthened version of the Wakamatsu lemma in the exact context, and also prove that subcategories which are $\mathcal{E}$-resolving and epimorphic precovering
with kernels in their right $\mathcal{E}$-orthogonal class and subcategories which are $\mathcal{E}$-coresolving and monomorphic preenveloping with cokernels
in their left $\mathcal{E}$-orthogonal class are determined by each other. As applications, we get some complete hereditary $\mathcal{E}$-cotorsion pairs.

Let $R$ be an associative ring with identity, and let $\mathscr{A}=\Mod R$ be the category of right $R$-modules and $\mathcal{E}$
contain all pure short exact sequences. As applications of results obtained in Section 3, we prove in Section 4 that for any subcategory $\mathcal{X}$ of $\Mod R$ in which
all modules are pure projective (resp. pure injective), the right
(resp. left) orthogonal class with respect to the exact structure $\mathcal{E}$ of $\mathcal{X}$ is preenveloping (resp. covering). Moreover,
we construct some complete hereditary $\mathcal{E}$-cotorsion pairs induced by the category of pure injective modules. Some results of G\"obel and Trlifaj in [11]
are obtained as corollaries.

\vspace{0.5cm}

{\bf 2. Exact Categories, Balanced Pairs and Relative Homology}

\vspace{0.2cm}

The following two definitions are cited from [4], see also [17] and [14].

\vspace{0.2cm}

{\bf Definition 2.1.} Let $\mathscr{B}$ be an additive category. A {\it kernel-cokernel pair} $(i,p)$ in $\mathscr{B}$ is a pair
of composable morphisms $$L\buildrel {i}\over \longrightarrow M\buildrel {p} \over \longrightarrow N$$ such that
$i$ is a kernel of $p$ and $p$ is a cokernel of $i$. If a class $\mathcal{E}$ of kernel-cokernel pairs on $\mathscr{B}$
is fixed, an {\it  admissible monic} (sometimes called {\it inflation}) is a morphism $i$ for which there exists a morphism $p$
such that $(i,p)\in \mathcal{E}$. {\it Admissible epics} (sometimes called {\it deflations}) are defined dually.

An \emph{exact structure} on $\mathscr{B}$ is a class $\mathcal{E}$ of kernel-cokernel pairs which is closed under isomorphisms and satisfies the following axioms:

[E0] For any object $B$ in $\mathscr{B}$, the identity morphism $\textrm{id}_B$ is both an admissible monic and an admissible epic.

[E1] The class of admissible monics is closed under compositions.

[{E1}$^{op}$] The class of admissible epics is closed under compositions.

[E2] The push-out of an admissible monic along an arbitrary morphism exists and yields an admissible monic, that is, for any admissible
monic $i:L\to M$ and any morphism $f:L\to M'$, there is a push-out diagram
$$\xymatrix{L\ar[r]^i\ar[d]_f & M\ar[d]^{f'}\\
M'\ar[r]^{i'}& X
}$$
with $i'$ an admissible monic.

[${\text E2^{op}}$] The pull-back of an admissible epic along an arbitrary morphism exists and yields an admissible epic, that is,
for any admissible epic $p:M\to N$ and any morphism $g:M'\to N$, there is a pull-back diagram
$$\xymatrix{X\ar[r]^{p'}\ar[d]_{G'} & M'\ar[d]^{g}\\
M\ar[r]^{p}& N
}$$
with $p'$ an admissible epic.

An {\it exact category} is a pair $(\mathscr{B},\mathcal{E})$ consisting of an additive category $\mathscr{B}$ and an exact structure $\mathcal{E}$ on $\mathscr{B}$.
Elements of $\mathcal{E}$ are called {\it admissible short exact sequences} (or {\it conflations}).

\vspace{0.2cm}

{\bf Definition 2.2.} Let $(\mathscr{B},\mathcal{E})$ be an exact category.

(1) An object $P\in\mathscr{B}$ is an ($\mathcal{E}$-){\it projective object} if for any admissible epic $p:M\to N$ and any morphism $f:P\to N$
there exists $f':P\to M$ such that $f=pf'$; an object $I\in\mathscr{B}$ is an ($\mathcal{E}$-){\it injective object} if for any admissible monic
$i:L\to M$ and any morphism $g:L\to I$ there exists $g':M\to I$ such that $g=g'i$.

(2) $(\mathscr{B},\mathcal{E})$ is said to {\it have enough projective objects} if for any object $M\in\mathscr{B}$ there exists an admissible epic
$p:P\to M$ with $P$ a projective object of $\mathscr{B}$; ($\mathscr{B},\mathcal{E}$) is said to {\it have enough injective objects} if for any
object $M\in\mathscr{B}$ there exists an admissible monic $i:M\to I$ with $I$ an injective object of $\mathscr{B}$.

\vspace{0.2cm}

We have the following standard observation.

\vspace{0.2cm}

{\bf Lemma 2.3.} {\it Let $(\mathscr{B},\mathcal{E})$ be an exact category with enough projective objects and enough injective objects,
and let $$0\to X \to Y\to Z\to 0\eqno{(2.1)}$$ be a sequence of morphisms in $\mathscr{B}$. Then the following statements are equivalent.

(1) $(2.1)$ is a conflation.

(2) For any projective object $P$ of $\mathscr{B}$, the induced sequence of abelian groups
$$0\to \Hom_{\mathscr{B}}(P,X)\to \Hom_{\mathscr{B}}(P,Y)\to \Hom_{\mathscr{B}}(P,Z)\to 0$$ is exact.

(3) For any injective object $I$ of $\mathscr{B}$, the induced sequence of abelian groups
$$0\to \Hom_{\mathscr{B}}(Z,I)\to \Hom_{\mathscr{B}}(Y,I)\to \Hom_{\mathscr{B}}(X,I)\to 0$$ is exact.}

\vspace{0.2cm}

{\it Proof.} The implications $(1)\Rightarrow(2)$ and $(1)\Rightarrow(3)$ are trivial. We just need to prove the implication
$(2)\Rightarrow(1)$ since the implication $(3)\Rightarrow(1)$ is its dual.

Let $$0\to X\buildrel {f_{2}} \over\to Y\buildrel {f_{1}} \over\to Z\to 0$$ be a sequence in
$\mathscr{B}$. One easily proves that $f_{2}$ is a monomorphism and that $f_{1}f_{2}=0$. We claim that $f_{2}$ is a kernel of $f_{1}$. In fact, for any $K\in \mathscr{B}$,
there exists an exact sequence $$P_{1}\to P_{0}\to K\to 0$$ with $P_{0}, P_{1}$ projective in $\mathscr{B}$. By (2) and the snake lemma, we obtain the following
commutative diagram with exact columns and rows:
$$\xymatrix{& 0\ar[d] & 0\ar[d] & 0\ar[d] & \\
0\ar@{-->}[r] & \Hom_{\mathscr{B}}(K,X)\ar@{-->}[r]^{\Hom_{\mathscr{B}}(K,f_{2})}\ar[d] & \Hom_{\mathscr{B}}(K,Y)
\ar@{-->}[r]^{\Hom_{\mathscr{B}}(K,f_{1})}\ar[d] & \Hom_{\mathscr{B}}(K,Z)\ar[d] & \\
0\ar[r] & \Hom_{\mathscr{B}}(P_0,X)\ar[r]^{\Hom_{\mathscr{B}}(P_0,f_{2})}\ar[d] & \Hom_{\mathscr{B}}(P_{0},Y)
\ar[r]^{\Hom_{\mathscr{B}}(P_0,f_{1})}\ar[d] & \Hom_{\mathscr{B}}(P_{0},Z)\ar[d]\ar[r] & 0 \\
0\ar[r] & \Hom_{\mathscr{B}}(P_1,X)\ar[r]^{\Hom_{\mathscr{B}}(P_1,f_{2})} & \Hom_{\mathscr{B}}(P_{1},Y)
\ar[r]^{\Hom_{\mathscr{B}}(P_1,f_{1})} & \Hom_{\mathscr{B}}(P_{1},Z)\ar[r] & 0.}$$
For any $u: K\to Y$, if $f_{1}u=0$, then
$$\Hom_{\mathscr{B}}(K,f_1)(u)=\Hom_{\mathscr{B}}(K,f_1u)=0$$ and
$u\in \Ker \Hom_{\mathscr{B}}(K,f_{1})=\Im \Hom_{\mathscr{B}}(K,f_{2})$. So there exists a unique morphism $v: K\to X$ such that
$u=\Hom_{\mathscr{B}}(K,f_{2})(v)= f_{2}v$. Thus $f_{2}$ is a kernel of $f_{1}$.

Again by (2), there exists a deflation $\varphi: P_{2}\to Z$ with $P_{2}$ projective in $\mathscr{B}$ such that $\varphi= f_{1}s$ for some
$s: P_{2}\to Y$. Then $f_{1}$ is a deflation by [4, Proposition 2.16]. So $$0\to X \to Y\to Z\to 0$$ is a conflation, as required.
\hfill{$\square$}

\vspace{0.2cm}

{\bf Definition 2.4.} ([1]) Let $\mathscr{A}$ be an abelian category with small Ext groups (that is, $\Ext^{1}_{\mathscr{A}}(X,Y)$ is a set for any $X, Y\in \mathscr{A}$).
For a subfunctor $\mathcal{F}\subseteq \Ext^{1}_{\mathscr{A}}(-,-)$, an object $C$ (resp. $D$) in $\mathscr{A}$ is called {\it $\mathcal{F}$-projective}
(resp. {\it $\mathcal{F}$-injective}) if $\mathcal{F}(C,-)=0$ (resp. $\mathcal{F}(-,D)=0$). An exact sequence $$0\to A' \to A \to A''\to 0$$ in $\mathscr{A}$
is called an {\it $\mathcal{F}$-sequence} if it is an element of $\mathcal{F}(A'',A')$.
The category $\mathscr{A}$ is said to {\it have enough $\mathcal{F}$-projective objects} if for any
$A\in \mathscr{A}$, there exists an $\mathcal{F}$-sequence
$$0\to X \to C\to A\to 0$$ such that
$C$ is $\mathcal{F}$-projective; and $\mathscr{A}$ is said to {\it have enough $\mathcal{F}$-injective objects} if for any $A\in \mathscr{A}$, there exists an $\mathcal{F}$-sequence
$$0\to A \to D\to Y\to 0$$ such that $D$ is $\mathcal{F}$-injective.

\vspace{0.2cm}

{\bf Definition 2.5.} ([6]) Let $\mathscr{C}$ be a subcategory of $\mathscr{A}$. A morphism $f: C\rightarrow D$ in $\mathscr{A}$
with $C \in\mathscr{C}$ is called a {\it $\mathscr{C}$-precover} of $D$ if for any morphism $g: C'\rightarrow D$ in $\mathscr{A}$
with $C' \in\mathscr{C}$, there exists a morphism $h:C'\rightarrow C$ such that $g=fh$.
The morphism $f:C\rightarrow D$ is called {\it right minimal} if an endomorphism $h:C\rightarrow C$ is an automorphism
whenever $f=fh$. A {\it $\mathscr{C}$-precover} is called a {\it $\mathscr{C}$-cover} if it is right minimal;
$\mathscr{C}$ is called a {\it (pre)covering subcategory} of $\mathscr{A}$ if every object in $\mathscr{A}$ has a $\mathscr{C}$-(pre)cover;
$\mathscr{C}$ is called an {\it epimorphic (pre)covering subcategory} of $\mathscr{A}$ if every object in $\mathscr{A}$ has an epimorphic
$\mathscr{C}$-(pre)cover. Dually, the notions of a {\it $\mathscr{C}$-(pre)envelope}, a {\it (pre)enveloping subcategory} and
a {\it monomorphic (pre)enveloping subcategory} are defined.

\vspace{0.2cm}

Let $\mathscr{C}$ be a subcategory of $\mathscr{A}$. Recall that a sequence in $\mathscr{A}$ is called {\it $\Hom_{\mathscr{A}}(\mathscr{C},-)$-exact}
if it is exact after applying the functor $\Hom_{\mathscr{A}}(C,-)$ for any object $C\in \mathscr{C}$. Let $M\in\mathscr{A}$.
An exact sequence (of finite or infinite length)
$$\cdots\buildrel {f_{n+1}} \over \longrightarrow C_{n}\buildrel {f_{n}} \over \longrightarrow\cdots \buildrel {f_{2}}
\over \longrightarrow C_{1}\buildrel {f_{1}} \over \longrightarrow C_{0}\buildrel {f_{0}} \over \longrightarrow M\to 0$$
in $\mathscr{A}$ with all $C_{i}\in\mathscr{C}$ is called a $\mathscr{C}${\it-resolution} of $M$ if it is $\Hom_{\mathscr{A}}(\mathscr{C},-)$-exact,
that is, each $f_{i}$ is an epimorphic $\mathscr{C}$-precover
of $\Im f_{i}$. We denote sometimes the $\mathscr{C}$-resolution of $M$ by $\mathscr{C}^{\bullet}\to M$,
where $$\mathscr{C}^{\bullet}:=\cdots \rightarrow C_{2}\buildrel {f_{2}} \over \longrightarrow C_{1}\buildrel {f_{1}} \over \longrightarrow C_{0}\to 0$$
is the {\it deleted $\mathscr{C}$-resolution} of $M$. Note that by a version of the comparison theorem, the $\mathscr{C}$-resolution is unique up to
homotopy ([7, p.169]). Dually, the notions of a $\Hom_{\mathscr{A}}(-, \mathscr{C})$-exact sequence and a $\mathscr{C}$-coresolution are defined.

\vspace{0.2cm}

{\bf Definition 2.6.} ([5]) A pair $(\mathscr{C},\mathscr{D})$ of additive subcategories in $\mathscr{A}$
is called a {\it balanced pair} if the following conditions are satisfied.

(1) $\mathscr{C}$ is epimorphic precovering and $\mathscr{D}$ is monomorphic preenveloping.

(2) For any $M\in \mathscr{A}$, there is a $\mathscr{C}$-resolution $\mathscr{C}^{\bullet}\to M$ such that it is $\Hom_{\mathscr{A}}(-,\mathscr{D})$-exact.

(3) For any $N\in \mathscr{A}$, there is a $\mathscr{D}$-coresolution $N\to \mathscr{D}^{\bullet}$ such that it is $\Hom_{\mathscr{A}}(\mathscr{C},-)$-exact.

\vspace{0.2cm}

{\bf Remark 2.7.} By [5, Proposition 2.2], in Definition 2.6, keeping condition (1) unaltered, the conditions (2) and (3) can be replaced by a common
condition $(2'+3')$: A short exact sequence $$0\to Y \to Z\to X\to 0$$ in $\mathscr{A}$ is $\Hom_{\mathscr{A}}(\mathscr{C},-)$-exact if and only if it
is $\Hom_{\mathscr{A}}(-,\mathscr{D})$-exact.

\vspace{0.2cm}

Some examples of balanced pairs are as follows.

\vspace{0.2cm}

{\bf Example 2.8.}

Let $R$ be an associative ring with identity and $\Mod R$ the category of right $R$-modules.

(1) ($\mathcal{P}_{0},\mathcal{I}_{0}$) is the standard balanced pair in $\Mod R$, where $\mathcal{P}_{0}$ and
$\mathcal{I}_{0}$ are the subcategories of $\Mod R$ consisting of projective modules and injective modules respectively.

(2) ($\mathcal{PP},\mathcal{PI}$) is a balanced pair in $\Mod R$ ([7, Example 8.3.2]), where $\mathcal{PP}$ and
$\mathcal{PI}$ are the subcategories of $\Mod R$ consisting of pure projective modules and pure injective modules respectively.

(3) If $R$ is a Gorenstein ring, then ($\mathcal{GP},\mathcal{GI}$) is a balanced pair in $\Mod R$ ([7, Theorem 12.1.4]), where $\mathcal{GP}$ and
$\mathcal{GI}$ are the subcategories of $\Mod R$ consisting of Gorenstein projective modules and Gorenstein injective modules respectively.

(4) Let $\mathscr{A}$ be an abelian category with enough projective and injective objects. If both ($\mathcal{B},\mathcal{C}$) and ($\mathcal{C},\mathcal{D}$)
are complete hereditary classical cotorsion pairs in $\mathscr{A}$,
then ($\mathcal{B},\mathcal{D}$) is a balanced pair in $\mathscr{A}$ ([5, Proposition 2.6]).

\vspace{0.2cm}

The main result in this section is the following

\vspace{0.2cm}

{\bf Theorem 2.9.} {\it Let $\mathscr{A}$ be an abelian category with small $\Ext$ groups. Then there exists a one-to-one correspondence between any two of the following.

(1) Balanced pairs $(\mathscr{C},\mathscr{D})$ in $\mathscr{A}$.

(2) Subfunctors $\mathcal{F}\subseteq \Ext^{1}_{\mathscr{A}}(-,-)$ such that $\mathscr{A}$ has enough $\mathcal{F}$-projectives and enough $\mathcal{F}$-injectives.

(3) Quillen exact structures $\mathcal{E}$ in $\mathscr{A}$ with enough $\mathcal{E}$-projectives and enough $\mathcal{E}$-injectives (that is, such that the resulting exact
category $(\mathscr{A},\mathcal{E})$ has enough projective and enough injective objects).}

\vspace{0.2cm}

{\it Proof.} We first show that there exists a one-to-one correspondence between (1) and (2).
Let ($\mathscr{C},\mathscr{D}$) be a balanced pair in $\mathscr{A}$. If we define $\mathcal{F}(X,Y)$ as the subset of $\Ext^{1}_{\mathscr{A}}(X,Y)$
consisting of the equivalence classes of short exact sequences of condition $(2'+3')$ in Remark 2.7, then the assignment $(X,Y)\longmapsto \mathcal{F}(X,Y)$
is the definition on objects of a subfunctor of $\Ext^{1}_{\mathscr{A}}(-,-)$, and \{$\mathcal{F}$-projectives\} = $\mathscr{C}$, \{$\mathcal{F}$-injectives\} = $\mathscr{D}$.
The converse is easy if we put ($\mathscr{C},\mathscr{D}$)= (\{$\mathcal{F}$-projectives\}, \{$\mathcal{F}$-injectives\}).

Next we show that there exists a one-to-one correspondence between (2) and (3).
Let $\mathcal{F}$ be a subfunctor of $\Ext^{1}_{\mathscr{A}}(-,-)$ such that $\mathscr{A}$ has enough $\mathcal{F}$-projectives and enough $\mathcal{F}$-injectives.
Put $\mathcal{E}= \{\mathcal{F}$-sequences\}, and \{$\mathcal{E}$-projectives\} = \{$\mathcal{F}$-projectives\}, \{$\mathcal{E}$-injectives\} = \{$\mathcal{F}$-injectives\}.
We claim that ($\mathscr{A},\mathcal{E}$) is an exact category with enough projective and enough injective objects.

Clearly, $\mathcal{E}$ contains \{$0\to A\buildrel {1_{A}}\over\to A\to 0 \to 0$\} and \{$0\to 0\to A\buildrel {1_{A}}\over\to A\to 0$\} for any
$A\in \mathscr{A}$.

For [E1], let $$0\to A_{1}\buildrel{f}\over\to A_{2}\to L\to 0\ {\rm and}$$ $$0\to A_{2}\buildrel{g}\over\to A_{3}\to K\to 0$$ be $\mathcal{F}$-sequences
in $\mathscr{A}$. We will prove $$0\to A_{1}\buildrel{gf}\over\to A_{3}\to \Coker gf\to 0$$ is an $\mathcal{F}$-sequence. In fact,
for any $\mathcal{F}$-injective object $I$ of $\mathscr{A}$ and $\alpha: A_{1}\to I$, by [1, Proposition 1.5] there exists $\beta: A_{2}\to I$ such that $\alpha = \beta f$.
Consider the following commutative diagram:
$$\xymatrix{& & 0\ar[d] & & & & \\
0\ar[r] & A_{1}\ar[r]^{f}\ar@{=}[d] & A_{2}\ar[r]\ar[d]^{g} & L\ar[r]\ar[d] & 0 & \\
0\ar[r] & A_{1}\ar[r]^{gf}\ar[d]^{\alpha} & A_{3}\ar[r]\ar[d]\ar@{-->}[ld]^{\theta}& \Coker gf\ar[r] &  0 &\\
& I & K\ar[d] & & & &\\
& &0. & & & &}$$
Then there exists $\theta: A_{3}\to I$ such that $\beta = \theta g$, that is, $\alpha = \beta f = \theta (gf)$. So $$0\to A_{1}\buildrel{gf}\over\to A_{3}\to \Coker gf\to 0$$
is an $\mathcal{F}$-sequence by [1, Proposition 1.5].

For [E2], take an $\mathcal{F}$-sequence $$0\to L\buildrel{i}\over\to M\to N\to 0.$$ For any morphism $f: L\to M'$, consider the following pushout diagram:
$$\xymatrix{0\ar[r] & L\ar[d]_{f}\ar[r]^{i} & M\ar[d]\ar[r] & N \ar@{=}[d]\ar[r] & 0 & \\
0\ar[r] & M'\ar[r]^{i'} & X\ar[r] & N\ar[r] & 0.}$$
By [1, Proposition 1.5], for any $\mathcal{F}$-injective object $I$ of $\mathscr{A}$ and $s: M'\to I$, there exists $t: M\to I$ such that $sf= ti$. It follows from the universal property of pushouts,
there exists $h: X\to I$ such that $s= hi'$. So $$0\to M'\buildrel{i'}\over\to X\to N\to 0$$ is an $\mathcal{F}$-sequence by [1, Proposition 1.5].
Dually [E1$^{op}$] and [E2$^{op}$] follow.

Conversely, for any $X, Y\in \mathscr{A}$, we define $\mathcal{F}(X,Y)$ as the subset of $\Ext^{1}_{\mathscr{A}}(X,Y)$ consisting of
the equivalence classes of admissible short exact sequences. Let $X',Y'\in \mathscr{A}$, $f: Y\to Y'$ and $g: X'\to X$. Then by [E2] and [E2$^{op}$], we get the following two commutative diagrams:
$$\xymatrix{\mathcal{F}(X,Y)\ar[r]^{\rm inc.}\ar[d]^{\mathcal{F}(X,f)} & \Ext_{\mathscr{A}}^{1}(X,Y)\ar[d]^{\Ext_{\mathscr{A}}^{1}(X,f)}\\
\mathcal{F}(X,Y')\ar[r]^{\rm inc.} & \Ext_{\mathscr{A}}^{1}(X,Y')},\ \mathrm{and}\ \ \
\xymatrix{\mathcal{F}(X,Y)\ar[r]^{\rm inc.}\ar[d]^{\mathcal{F}(g,Y)} & \Ext_{\mathscr{A}}^{1}(X,Y)\ar[d]^{\Ext_{\mathscr{A}}^{1}(g,Y)}\\
\mathcal{F}(X',Y)\ar[r]^{\rm inc.} & \Ext_{\mathscr{A}}^{1}(X',Y).}$$
So $\mathcal{F}$ is a subfunctor of $\Ext^{1}_{\mathscr{A}}(-,-)$
such that $\mathscr{A}$ has enough $\mathcal{F}$-projectives and enough $\mathcal{F}$-injectives. \hfill{$\square$}

\vspace{0.5cm}

{\bf 3. $\mathcal{E}$-Cotorsion Pairs}

\vspace{0.2cm}

In this section, a pair ($\mathscr{A},\mathcal{E}$) means that $\mathscr{A}$ is an abelian category with small
Ext groups and $\mathcal{E}$ is an exact structure on $\mathscr{A}$ such that $\mathscr{A}$ has enough $\mathcal{E}$-projectives
and enough $\mathcal{E}$-injectives. We first give a strengthened version of the Wakamatsu lemma in the exact context,
and then apply it to obtain complete hereditary $\mathcal{E}$-cotorsion pairs in ($\mathscr{A},\mathcal{E}$).
For any $M, N\in \mathscr{A}$ and $i\geq 1$,
we use $\mathrm{\mathcal{E}xt}^{i}_{\mathscr{A}}(M,N)$ to denote the $i$-th cohomology group by taking $\mathcal{E}$-projective
resolution of $M$ or taking $\mathcal{E}$-injective coresolution of $N$.

Inspired by [20], we give the following

\vspace{0.2cm}

{\bf Definition 3.1.} Let ($\mathscr{A},\mathcal{E}$) be a pair as above.

(1) A pair $(\mathscr{X},\mathscr{Y})$ of full subcategories of $\mathscr{A}$ is called an {\it $\mathcal{E}$-cotorsion pair}
provided that $\mathscr{X}={^{\bot_{\ast}}\mathscr{Y}}$ and $\mathscr{Y}={\mathscr{X}^{\bot_{\ast}}}$,
where $${\mathscr{X}^{\bot_{\ast}}}=\{N\in \mathscr{A} \mid  \mathrm{\mathcal{E}xt}^{1}_{\mathscr{A}}(M,N)=0\  \mathrm{for}\  \mathrm{any}\   M\in \mathscr{X}\},\ {\rm and}$$
$${^{\bot_{\ast}}\mathscr{Y}}=\{M\in \mathscr{A} \mid  \mathrm{\mathcal{E}xt}^{1}_{\mathscr{A}}(M,N)=0\  \mathrm{for}\ \mathrm{any}\  N\in\mathscr{Y}\}.$$

(2) An $\mathcal{E}$-cotorsion pair $(\mathscr{X},\mathscr{Y})$ is called {\it complete} if for any $M\in \mathscr{A}$, there exist two conflations of the form:
$$0\to Y\to X\to M\to 0,\ {\rm and}$$
$$0\to M\to Y'\to X'\to 0 $$ with $X, X'\in \mathscr{X}$ and $Y, Y'\in \mathscr{Y}$.

(3) A subcategory $\mathscr{T}$ of $\mathscr{A}$ is called {\it closed under $\mathcal{E}$-extensions} if the end terms in a conflation are in $\mathscr{T}$
implies that the middle term is also in $\mathscr{T}$. The subcategory $\mathscr{T}$
is called {\it $\mathcal{E}$-resolving} if it contains all $\mathcal{E}$-projectives of $\mathscr{A}$, closed under $\mathcal{E}$-extensions, 
and for any conflation $$0\to X\to Y\to Z\to 0$$ with $Y, Z\in \mathscr{T}$, we have $X\in \mathscr{T}$.
Dually, the notion of {\it $\mathcal{E}$-coresolving subcategories} is defined.

(4) An $\mathcal{E}$-cotorsion pair ($\mathscr{X},\mathscr{Y}$) is called {\it hereditary} if $\mathscr{X}$ is $\mathcal{E}$-resolving.

\vspace{0.2cm}

{\it Remark 3.2.} (1) By [11, Lemma 2.2.10], an $\mathcal{E}$-cotorsion pair ($\mathscr{X},\mathscr{Y}$) is hereditary if and only if $\mathscr{Y}$ is $\mathcal{E}$-coresolving.

(2) If $\mathcal{E}$ is the abelian structure, then all the notions in Definition 3.1 coincide with the classical ones in [5, p.6].

(3) \{$\mathcal{E}$-projectives\} = ${^{\bot_{\ast}}\mathscr{A}}$ and \{$\mathcal{E}$-injectives\} = ${\mathscr{A}^{\bot_{\ast}}}$.

\vspace{0.2cm}

The following result is a strengthened version of the translation of the Wakamatsu lemma to the exact context.

\vspace{0.2cm}

{\bf Theorem 3.3.} {\it Let $(\mathscr{A},\mathcal{E})$ be a pair as above, $\mathscr{X}$ a full subcategory of $\mathscr{A}$ closed under $\mathcal{E}$-extensions and $A\in \mathscr{A}$.

(1) If $f$: $X\rightarrow A$ is an epimorphic $\mathscr{X}$-cover, then $\Ker f\in \mathscr{X}^{\bot_{\ast}}$.

(2) If $f$: $A\rightarrow X$ is a monomorphic $\mathscr{X}$-envelope, then $\Coker f\in {^{\bot_{\ast}}\mathscr{X}}$.}

\vspace{0.2cm}

{\it Proof.} We only prove (2), and (1) is dual to (2).

By assumption, there exists an exact sequence $$0\to A\buildrel {f} \over \rightarrow X\buildrel {g} \over \rightarrow \Coker f \to 0$$ in $\mathscr{A}$
with $f:A\to X$ a monomorphic $\mathscr{X}$-envelope. To prove $\Coker f \in {^{\bot_*}\mathscr{X}}$,
it suffices to prove that for any $X'\in\mathscr{X}$, any conflation $$0\to X'\to M\buildrel {h} \over \rightarrow \Coker f \to 0$$ in $\mathscr{A}$ splits.
Consider the pullback of $g$ and $h$:
$$\xymatrix{& & 0\ar[d] & 0\ar[d] & & & \\
& & X'\ar@{=}[r]\ar[d] & X'\ar[d] & & & \\
0\ar[r] & A\ar[r]\ar@{=}[d] & Y\ar[r]\ar[d] & M\ar[r]\ar[d]^{h} & 0 & \\
0\ar[r] & A\ar[r]^{f} & X\ar[r]^{g}\ar[d] & \Coker f \ar[r]\ar[d] &  0 &\\
& &0 & 0. & & &}$$
Then the middle column is a conflation. Since $X', X\in\mathscr{X}$ and $\mathscr{X}$ is closed under $\mathcal{E}$-extensions, we have $Y\in\mathscr{X}$.
Because $f: A\to X$ is an $\mathscr{X}$-envelope of $A$, we obtain the following commutative diagram with exact rows:
$$\xymatrix{0\ar[r] & A\ar@{=}[d]\ar[r]^{f} & X\ar[d]\ar[r] & \Coker f \ar[d]\ar[r] & 0 & \\
0\ar[r] & A\ar@{=}[d]\ar[r] & Y\ar[d]\ar[r] & M\ar[d]^{h}\ar[r] & 0 & \\
0\ar[r] & A\ar[r]^{f} & X\ar[r] & \Coker f \ar[r] & 0. &}$$
Since $f$ is left minimal, we have that the composition $$X\to Y\to X$$ is an isomorphism.
This implies that the composition $$\Coker f \to M\buildrel {h} \over \rightarrow \Coker f$$ is also an isomorphism.
Therefore the conflation $$0\to X'\to M\buildrel {h} \over \rightarrow \Coker f \to 0$$ splits and the assertion follows. \hfill{$\square$}

\vspace{0.2cm}

{\bf Theorem 3.4.} {\it Let $(\mathscr{A},\mathcal{E})$ be a pair as above and $(\mathscr{X},\mathscr{Y})$ a pair of full subcategories of $\mathscr{A}$.
Then the following statements are equivalent.

(1) $\mathscr{X}$ is $\mathcal{E}$-resolving, $\mathscr{Y}=\mathscr{X}^{\bot_{\ast}}$ and each $A\in \mathscr{A}$
admits an epimorphic $\mathscr{X}$-precover $p$: $X\rightarrow A$ such that $\Ker p\in \mathscr{Y}$.

(2) $\mathscr{Y}$ is $\mathcal{E}$-coresolving, $\mathscr{X}={^{\bot_{\ast}}\mathscr{Y}}$ and each $A\in \mathscr{A}$
admits a monomorphic $\mathscr{Y}$-preenvelope $j$: $A\rightarrow Y$ such that $\Coker j\in \mathscr{X}$.

(3) $(\mathscr{X},\mathscr{Y})$ is a complete hereditary $\mathcal{E}$-cotorsion pair.}

\vspace{0.2cm}

{\it Proof.} $(2)\Rightarrow(1)$ First observe that $\mathscr{X}(={^{\bot_{\ast}}\mathscr{Y}})$ is closed under direct summands,
$\mathcal{E}$-extensions and contains all $\mathcal{E}$-projective objects of $\mathscr{A}$. For any conflation
$$0\to A\to B\to C\to 0$$ in $\mathscr{A}$ with $B,C\in \mathscr{X}$, we have an exact sequence
$$0=\mathrm{\mathcal{E}xt}^{1}_{\mathscr{A}}(B,Y)\to \mathrm{\mathcal{E}xt}^{1}_{\mathscr{A}}(A,Y)\to \mathrm{\mathcal{E}xt}^{2}_{\mathscr{A}}(C,Y)$$ for any $Y\in \mathscr{Y}$.
Take a conflation
$$0\to Y\to E\to L\to 0$$ in $\mathscr{A}$ with $E$ $\mathcal{E}$-injective. We have $L\in \mathscr{Y}$ since $\mathscr{Y}$ is $\mathcal{E}$-coresolving. Then we have an exact sequence:
$$0=\mathrm{\mathcal{E}xt}^{1}_{\mathscr{A}}(C,E)\to \mathrm{\mathcal{E}xt}^{1}_{\mathscr{A}}(C,L)\to \mathrm{\mathcal{E}xt}^{2}_{\mathscr{A}}(C,Y)\to \mathrm{\mathcal{E}xt}^{2}_{\mathscr{A}}(C,E)=0.$$
Notice that $\mathrm{\mathcal{E}xt}^{1}_{\mathscr{A}}(C,L)=0$, so $\mathrm{\mathcal{E}xt}^{2}_{\mathscr{A}}(C,Y)=0$, and hence $\mathrm{\mathcal{E}xt}^{1}_{\mathscr{A}}(A,Y)=0$,
that is, $A\in \mathscr{X}$. Consequently we conclude that $\mathscr{X}$ is $\mathcal{E}$-resolving.

Next for each $A\in \mathscr{A}$, there exists a conflation $$0 \to K \to P \to A \to 0$$ in $\mathscr{A}$ with $P$ $\mathcal{E}$-projective.
By (2) and Lemma 2.3, we have a conflation
$$0 \to K \buildrel {j} \over\to Y \to M \to 0$$ with $Y\in \mathscr{Y}$ and $M \in\mathscr{X}$.
Consider the following pushout diagram:
$$\xymatrix{
& & 0\ar[d] & 0\ar[d]  & & &\\ & 0\ar[r] &
K\ar[r]\ar[d]^{j} & P\ar[r]\ar[d] &
A\ar[r]\ar@{=}[d] & 0 & & &\\ & 0\ar[r] &
Y\ar[r]\ar[d] & X\ar[r]^{p}\ar[d] & A\ar[r] & 0 & & &\\
& & M \ar@{=}[r]\ar[d]&M\ar[d] & & &\\& &0 &0. & & &}$$
Since $\mathscr{X}$ is $\mathcal{E}$-resolving and $P\in \mathscr{X}$,
it follows from the middle column in the above diagram that $X\in \mathscr{X}$. Obviously, the middle row $$0\to Y \to X \buildrel {p} \over\to A \to 0$$ is a conflation and
$Y\in \mathscr{Y}$. Therefore we deduce that $p$: $X\rightarrow A$ is an epimorphic $\mathscr{X}$-precover, as required.

It remains to prove that $\mathscr{Y}=\mathscr{X}^{\bot_{\ast}}$. Observe that $\mathscr{Y} \subseteq (^{\bot_{\ast}}\mathscr{Y})^{\bot_{\ast}}=\mathscr{X}^{\bot_{\ast}}$.
By (2) and Lemma 2.3, for any $A'\in \mathscr{X}^{\bot_{\ast}}=(^{\bot_{\ast}}\mathscr{Y})^{\bot_{\ast}}$, there exists a conflation $$0 \to A' \to Y' \to X' \to 0\eqno{(3.1)}$$
with $Y'\in \mathscr{Y}$ and $X'\in \mathscr{X}={^{\bot_{\ast}}\mathscr{Y}}$. So (3.1) splits. Thus $A'\in \mathscr{Y}$ since $\mathscr{Y}$ is closed under direct summands.

Dually, we get $(1)\Rightarrow(2)$.

The implications $(1)+ (2)\Rightarrow (3)$ and $(1)\Leftarrow (3)\Rightarrow (2)$ are clear. \hfill{$\square$}

\vspace{0.2cm}

As direct consequences of Theorems 3.3 and 3.4, we have the following

\vspace{0.2cm}

{\bf Corollary 3.5.} {\it

(1) If $\mathscr{X}$ is an epi-covering $\mathcal{E}$-resolving subcategory of $\mathscr{A}$, then each $A\in \mathscr{A}$
admits a monomorphic $\mathscr{Y}$-preenvelope $j$: $A\rightarrow Y$ with $\Coker j\in \mathscr{X}$, where $\mathscr{Y}=\mathscr{X}^{\bot_{\ast}}$.

(2) If $\mathscr{Y}$ is a mono-enveloping $\mathcal{E}$-coresolving subcategory of $\mathscr{A}$, then each $A\in \mathscr{A}$ admits
an epimorphic $\mathscr{X}$-precover $p$: $X\rightarrow A$ with $\Ker p\in \mathscr{Y}$, where $\mathscr{X}={^{\bot_{\ast}}\mathscr{Y}}$.}

\vspace{0.2cm}

{\bf Corollary 3.6.} {\it Let $(\mathscr{X}$,$\mathscr{Y})$ be a complete hereditary (classical) cotorsion pair in $\mathscr{A}$. Then we have

(1) If $\mathscr{X}$ contains the $\mathcal{E}$-projective objects of $\mathscr{A}$, then $(\mathscr{X}, \mathscr{X}^{\bot_{\ast}})$ is a complete hereditary
$\mathcal{E}$-cotorsion pair and each $A\in \mathscr{A}$ admits a monomorphic $\mathscr{X}^{\bot_{\ast}}$-preenvelope $j$: $A\rightarrow Y$ such that $\Coker j\in \mathscr{X}$.

(2) If $\mathscr{Y}$ contains the $\mathcal{E}$-injective objects of $\mathscr{A}$, then $({^{\bot_{\ast}}\mathscr{Y}}, \mathscr{Y})$ is a complete hereditary
$\mathcal{E}$-cotorsion pair and each $A\in \mathscr{A}$ admits an epimorphic ${^{\bot_{\ast}}\mathscr{Y}}$-precover $p$: $X\rightarrow A$ such that $\Ker p\in \mathscr{Y}$.}

\vspace{0.5cm}

{\bf 4. Applications}

\vspace{0.2cm}

In this section,  we will apply the results obtained in Section 3 to the module category.

Let $R$ be an associative ring with identity and $\mod R$ the category of finitely presented right $R$-modules.
Recall that a short exact sequence $$\xi:0\to A\to B\to C\to 0$$ in $\Mod R$ is called {\it pure exact} if the induced sequence
$\Hom_{R}(F,\xi)$ is exact for any $F\in\mod R$. In this case $A$ is called a {\it pure submodule} of $B$ and $C$ is called a
{\it pure quotient module} of $B$. In addition, modules that are projective (resp. injective) relative to pure exact sequences
are called {\it pure projective} (resp. {\it pure injective}). The subcategory of $\Mod R$ consisting of pure projective
(resp. pure injective) modules is denoted by $\mathcal{P}\mathcal{P}$ (resp. $\mathcal{P}\mathcal{I}$).

\vspace{0.2cm}

{\bf Lemma 4.1.} ([18, Corollary 3.5(c)]). {\it Let $\mathscr{F}$ be a subcategory of $\Mod R$ closed under pure submodules.
Then $\mathscr{F}$ is preenveloping if and only if $\mathscr{F}$ is closed under direct products.}

\vspace{0.2cm}

{\bf Lemma 4.2.} ([12, Theorem 2.5]). {\it Let $\mathscr{F}$ be a subcategory of $\Mod R$ closed under pure quotient modules.
Then the following statements are equivalent.

(1) $\mathscr{F}$ is closed under direct sums.

(2) $\mathscr{F}$ is precovering.

(3) $\mathscr{F}$ is covering.}

\vspace{0.2cm}

{\bf Lemma 4.3.} {\it Let $(\mathscr{A},\mathcal{E})$ be a pair as in Section 3 and $A\in \mathscr{A}$. Then the following statements are equivalent.

(1) $\Ext_{\mathscr{A}}^{1}(-,A)$ vanishes on all $\mathcal{E}$-projective objects of $\mathscr{A}$.

(2) Any short exact sequence $$0\to A\to Z\to V\to 0$$ is a conflation.}

{\it Proof.} $(2)\Rightarrow(1)$ is trivial.

$(1)\Rightarrow(2)$ For any $\mathcal{E}$-projective object $P$ of $\mathscr{A}$ and $f: P\to V$, consider the pullback diagram:
$$\xymatrix{0\ar[r] & A\ar@{=}[d]\ar[r] & N\ar[d]^{\alpha}\ar[r]^{i} & P \ar[d]^{f}\ar[r] & 0 & \\
0\ar[r] & A\ar[r] & Z\ar[r]^{g} & V\ar[r] & 0. & }$$
By (1), there exists $j: P\to N$ such that $ij= 1_P$. Then we have $f= fij= g\alpha j$. So by Lemma 2.3, the short exact sequence
$$0\to A\to Z\to V\to 0$$ is a conflation.  \hfill{$\square$}

\vspace{0.2cm}

Recall from [3] that a subcategory of $\Mod R$ is called {\it definable} if it is closed under arbitrary direct products, direct limits,
and pure submodules.

\vspace{0.2cm}

{\bf Theorem 4.4.} {\it Let $\mathcal{E}$ be an exact structure on $\Mod R$ with enough $\mathcal{E}$-projectives
and enough $\mathcal{E}$-injectives such that $\mathcal{E}$ contains all pure short exact sequences. If $\mathscr{X}\subseteq \mathcal{P}\mathcal{P}$
and $\mathscr{Y}\subseteq \mathcal{P}\mathcal{I}$ are two subcategories, then we have the following

(1) $\mathscr{X}^{\bot_{\ast}}$ is preenveloping. If moreover, each $X\in \mathscr{X}$ admits an $\mathcal{E}$-presention
$$C_{2}\to C_{1}\to C_{0}\to X\to 0,$$ where the $C_{i}$ are $\mathcal{E}$-projective in $\Mod R$ and finitely presented, then $\mathscr{X}^{\bot_{\ast}}$ is also covering.

(2) ${^{\bot_{\ast}}\mathscr{Y}}$ is covering and closed under pure quotients.

(3) If $\mathscr{Y}$ is closed under taking $\mathcal{E}$-cosyzygies, then we have

\ \ \ \ (i) $(^{\bot_{\ast}} \mathscr{Y},(^{\bot_{\ast}} \mathscr{Y})^{\bot_{\ast}})$ is a complete hereditary $\mathcal{E}$-cotorsion pair.

\ \ \ \ (ii) If $\Ext^{1}_{R}(-,Y)$ vanishes on $\mathcal{E}$-projective objects of $\Mod R$ for any $Y\in \mathscr{Y}$, then
$(^{\bot_{1}} \mathscr{Y},(^{\bot_{1}} \mathscr{Y})^{\bot_{\ast}})$ is a complete hereditary $\mathcal{E}$-cotorsion pair.}

\vspace{0.2cm}

{\it Proof.} (1) First observe that $\mathscr{X}^{\bot_{\ast}}$ is closed under direct products. Let
$$0\to A\to B\to C\to 0$$ be a pure exact sequence in $\Mod R$ with $B\in \mathscr{X}^{\bot_{\ast}}$. Then it is a conflation by assumption.
For any $X\in \mathscr{X}$, we have the following exact sequence:
$$\mathrm{Hom}_{R}(X,B) \buildrel {f} \over \longrightarrow \mathrm{Hom}_{R}(X,C)\to \mathrm{\mathcal{E}xt}^{1}_{\mathscr{A}}(X,A)\to \mathrm{\mathcal{E}xt}^{1}_{\mathscr{A}}(X,B)=0.$$
Notice that $f$ is epic, so $\mathrm{\mathcal{E}xt}^{1}_{\mathscr{A}}(X,A)=0$ and $A\in \mathscr{X}^{\bot_{\ast}}$. Thus $\mathscr{X}^{\bot_{\ast}}$ is closed under pure submodules.
It follows from Lemma 4.1 that $\mathscr{X}^{\bot_{\ast}}$ is preenveloping.

It is easy to see that $\mathscr{X}^{\bot_{\ast}}$ is closed under direct limits by assumption; in particular $\mathscr{X}^{\bot_{\ast}}$ is closed under direct sums.
Then $\mathscr{X}^{\bot_{\ast}}$ is definable. So $\mathscr{X}^{\bot_{\ast}}$ is closed under pure quotient modules by [3, Proposition 4.3(3)], and hence
$\mathscr{X}^{\bot_{\ast}}$ is covering by Lemma 4.2.

(2) Obviously, $^{\bot_{\ast}} \mathscr{Y}$ is closed under direct sums. By Lemma 4.2, it suffices to show that $^{\bot_{\ast}} \mathscr{Y}$
is closed under pure quotient modules. Let $$0\to A\to B\to C\to 0$$ be a pure exact sequence in $\Mod R$ with $B\in {^{\bot_{\ast}}\mathscr{Y}}$.
Then it is a conflation by assumption. By using an argument similar to that in the proof of (1), we get $C\in {^{\bot_{\ast}}\mathscr{Y}}$.

(3) (i) Because ${^{\bot_{\ast}}\mathscr{Y}}$ contains all projective modules, we have that $^{\bot_{\ast}} \mathscr{Y}$ is epimorphic covering by (2).
We claim that $^{\bot_{\ast}} \mathscr{Y}$ is $\mathcal{E}$-resolving. This will complete the proof of (i) by Theorems 3.3 and 3.4.

Clearly, ${^{\bot_{\ast}}\mathscr{Y}}$ contains all $\mathcal{E}$-projective objects of $\Mod R$ and ${^{\bot_{\ast}}\mathscr{Y}}$ is closed under direct summands
and $\mathcal{E}$-extensions. Now take a conflation $$0\to M\to N\to L\to 0$$ in $\Mod R$ with $N,L\in {^{\bot_{\ast}} \mathscr{Y}}$. By assumption, we have
$\Ker \mathrm{\mathcal{E}xt}^{1}_{\mathscr{A}}(-,Y)\subseteq \Ker \mathrm{\mathcal{E}xt}^{2}_{\mathscr{A}}(-,Y)$ for any $Y\in \mathscr{Y}$.
Now by the dimension shifting it yields that $M\in {^{\bot_{\ast}} \mathscr{Y}}$. So $^{\bot_{\ast}} \mathscr{Y}$ is $\mathcal{E}$-resolving.
The claim follows.

(ii) By Lemma 4.3, we have $^{\bot_{1}} \mathscr{Y}= {^{\bot_{\ast}}\mathscr{Y}}$. So, as a particular case of (i), the assertion follows. \hfill{$\square$}
\vspace{0.2cm}

Recall that a classical cotorsion pair ($\mathscr{X},\mathscr{Y}$) in $\Mod R$ is {\it perfect} if $\mathscr{X}$ is covering and $\mathscr{Y}$ is enveloping.
As a consequence of Theorem 4.4, we have the following

\vspace{0.2cm}

{\bf Corollary 4.5.}  {\it If $\mathscr{X}\subseteq \mathcal{P}\mathcal{P}$ and $\mathscr{Y}\subseteq \mathcal{P}\mathcal{I}$ are two subcategories in $\Mod R$, then we have

(1) $\mathscr{X}^{\bot_{1}}$ is preenveloping. If moreover, $R$ is a right coherent ring and $\mathscr{X}$ is a subcategory of $\mod R$, then $\mathscr{X}^{\bot_{1}}$ is also covering.

(2) ${^{\bot_{1}}\mathscr{Y}}$ is covering and closed under pure quotients.

(3) ([11, Theorem 3.2.9]) $(^{\bot_{1}} \mathscr{Y},(^{\bot_{1}} \mathscr{Y})^{\bot_{1}})$ is a perfect hereditary cotorsion pair.}

\vspace{0.2cm}

{\it Proof.} The assertions (1) and (2) follow directly from Theorem 4.4, and the assertion (3) follows from Theorem 4.4 and [7, Theorem 7.2.6]. \hfill{$\square$}

\vspace{0.5cm}

{\bf Acknowledgements.} This research was partially supported by NSFC (Grant Nos. 11171142 and 11571164)
and a Project Funded by the Priority Academic Program Development of Jiangsu Higher Education
Institutions. The authors would like to express their sincere thanks to the referees for many considerable suggestions, which have greatly improved this paper.

\end{document}